\documentclass[12pt]{article}
\usepackage{epsf}
\usepackage{amsmath}
\usepackage{amssymb}

\newtheorem{preExample}{\bf Example}
\newenvironment{Example}{\begin{preExample}\hspace{-0.5em}. }{\end{preExample}}
\newcommand{\xmpl}[1]{\begin{Example}{#1}\end{Example}}

\newtheorem{preLemma}{\bf Lemma}
\newenvironment{Lemma}{\begin{preLemma}\hspace{-0.5em}. }{\end{preLemma}}
\newcommand{\lem}[1]{\begin{Lemma}{#1}\end{Lemma}}

\newtheorem{preTHeorem}{\bf Theorem}
\newenvironment{THeorem}{\begin{preTHeorem}\hspace{-0.5em}. }{\end{preTHeorem}}
\newcommand{\thm}[1]{\begin{THeorem}{#1}\end{THeorem}}

\newtheorem{preDEfinition}{\bf Definition}
\newenvironment{DEfinition}{\begin{preDEfinition}\hspace{-0.5em}. }{\end{preDEfinition}}
\newcommand{\dfn}[1]{\begin{DEfinition}{#1}\end{DEfinition}}

\newenvironment{proof}{{\bf proof.}}{$\blacksquare$}
\newcommand{\prf}[1]{\begin{proof}#1\end{proof}}

\DeclareMathOperator{\cut}{Cut}
\DeclareMathOperator{\spec}{Spec}
\DeclareMathOperator{\code}{Code}
\DeclareMathOperator{\del}{Del}

\def\os{\Big\{}
\def\cs{\Big\}}
\def\wh{\widehat}
\def\wt{\widetilde}
\def\svert{\;\vert\;}
\def\G{{\cal G}}
\def\R{{\mathbb R}}
\def\d{{\mathrm d}}
\newcommand{\abs}[1]{\left\vert #1\right\vert}
\author{Hamed Daneshpajouh \and Hamid Reza Daneshpajouh
\thanks{Department of Mathematics and Computer Science, University of Tehran. 
Email: {\tt h.r.daneshpajouh@khayam.ut.ac.ir}. }
\and Farzad Didehvar
\thanks{Department of Mathematics and Computer science, Amirkabir University 
of Technology.
Email: {\tt
didehvar@aut.ac.ir, fr@ipm.ir}. }
}
\title{A metric on the space of weighted graphs}
\begin{document}
\maketitle
\begin{abstract}
In this paper we offer a metric similar to graph edit distance which 
measures the distance between two (possibly infinite)weighted graphs with 
finite norm (we define the norm of a graph as the sum of absolute values of its edges).
The main result is the completeness of the space. Some other 
analytical properties of this space are also investigated.
\end{abstract}
\section{Introduction}
Many objects can be demonstrated with weighted graphs.
In any collection of objects of similar nature a way
to quantify the difference between objects may be desired
(For instance if we were to select the most similar objects to a given object
from a database).
In the theoretical side one common way is to develop a metric on 
the space of objects in demand.
One way to build a metric, is to define some operations that transform
the members of the space to one another, and assign a cost to each operation
then define the distance between two objects to be the minimum cost that  
must be payed to transform the first object to the second via a sequence of the defined 
operations. Such metrics sometimes are referred to as ``Edit distance".
Two examples of them are the ``Levenshtein edit distance" \cite{1} 
on strings and ``Graph edit distance" \cite{2} on the space of finite graphs.
This paper extends the Graph edit distance to the space of 
``countable weighted graphs with finite norm" and investigates some topological
properties of the space.
\section{Priliminaries and intuitive examples}
In this chapter, we introduce the concepts intuitively. The main
question here is ``given two graphs, how much they differ?". Based on
this question we could define different distances, we choose here
``Graph edit distance" and we generalize it to infinite graphs.
\xmpl{
Consider the following two graphs
\newpage
\begin{figure}[ht]
\begin{center}
\epsfbox{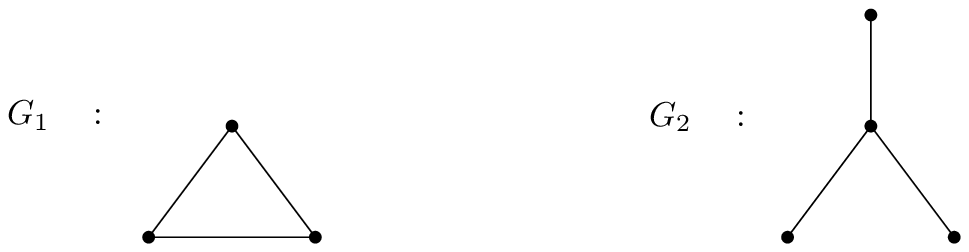}
\end{center}
\end{figure}
\noindent
One can transform $G_1$ to $G_2$, by adding a vertex and an edge
to $G_1$ and deleting an edge from it
\begin{figure}[ht]
\begin{center}
\epsfbox{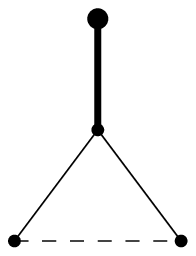}
\end{center}
\end{figure}
}
Given two graphs, it is possible to transform one to the other 
by addition and deletion of some vertices and edges. The minimum number of edge
addition and deletions in such a process is the distance between the two 
graphs and is denoted by $\d(G_1,G_2)$. In the above example $\d(G_1,G_2)=2$,
because we added an edge and deleted one.
It is clear that if two graphs differ only in 
isolated vertices, then by this definition their 
distance is zero.
\xmpl{
Consider the following weighted graphs
\begin{figure}[h]
\begin{center}
\epsfbox{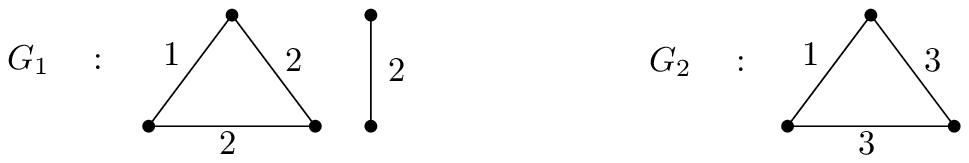}
\end{center}
\end{figure}\newline
We transform $G_1$ to $G_2$ as follows
\newpage
\begin{figure}[h]
\begin{center}
\epsfbox{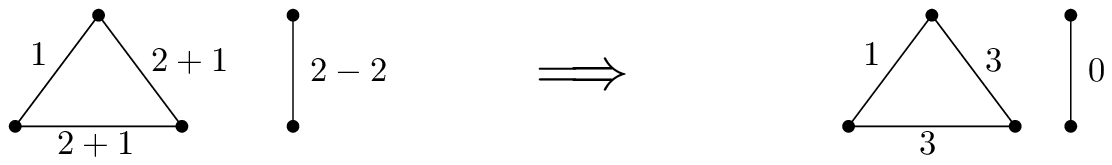}
\end{center}
\end{figure}
\noindent
The right-hand graph, after removing the zero weighted edge and it's endpoints
is same as $G_2$.
We define the distance between two edge-weighted graphs to be the minimum 
amount of edge-weight modifications required to transform one to the other
(zero-weighted edges and isolated vertices could be added and deleted for free)
, in this example $\d(G_1,G_2)=\abs{1}+\abs{1}+\abs{-2}=4$.
}
\xmpl{
We construct the sequence $\{G_n\}$ as follows
\begin{figure}[ht]
\begin{center}
\epsfbox{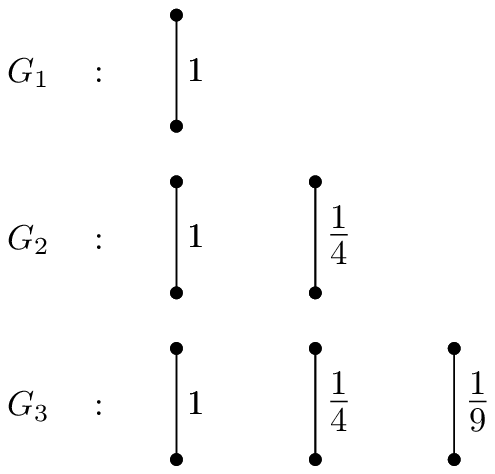}
\end{center}
\end{figure}
\newline
If $n<m$ then clearly
$$\d(G_n,G_m)\leq \sum_{t=n+1}^m\frac{1}{t^2}$$
hence $\{G_n\}$ is a Cauchy sequence. On the other hand $\{G_n\}$ does not 
approach to a finite graph (to prove that let $G$ be a graph with $m$ edges,
and show that for each $n>m$, $\d(G_n,G_m)\geq 1/(m+1)$), implying that the 
space of finite weighted graphs is not complete. But this sequence approaches 
to the following infinite graph
\begin{figure}[ht]
\begin{center}
\epsfbox{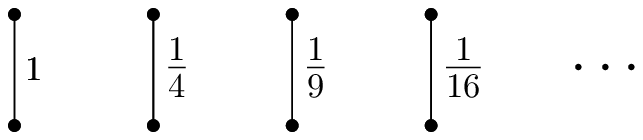}
\end{center}
\end{figure}
}
\xmpl{
Consider the following two graphs
\newpage
\begin{figure}[ht]
\begin{center}
\epsfbox{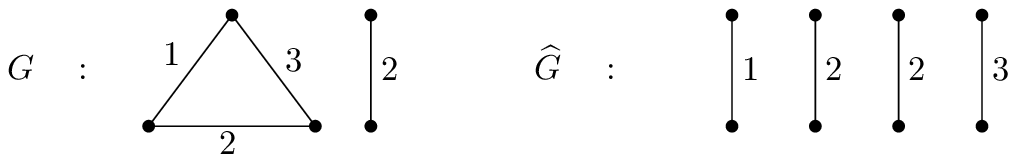}
\end{center}
\end{figure}
\noindent
$\wh G$ is obtained from $G$ by removing the joints.
You should be able to verify the following inequality for arbitrary graphs intuitively. 
$$\d(\wh{G},\wh{H})\leq \d(G,H)$$
}
\section{Developing the metric mathematically}
We first define the distance between two labeled graphs and then
define unlabeled graphs as the equivalence classes of labeled
graphs. The main result in this chapter is to show that the
introduced distance provides us a metric space.
\dfn{
Fix $V=\{ v_1,v_2,v_3,\cdots \}$ as the vertex set. We show the
set of all 2-element subsets of $V$ by $E$ (edge set).
A labeled graph is a function $w:E\to \R$. 
The value $w(e)$, is the weight of $e$(in $w$). 
A nonzero edge (of $w$), is one which it's weight is nonzero.
An isolated vertex of $w$ is a vertex that has no nonzero edge attached. 
$w$ is standard if it has infinitely many isolated vertices. 
We denote by $W$ the set of all labeled graphs, and by $W^s$ the set of 
all standard labeled graphs. 
For $w,w'\in W$, the norm of $w$, the distance between $w$ and $w'$ and 
the sets $W_0$ and $W^s_0$ are defined as follows
\begin{eqnarray*}\begin{array}{c}
|w|=\sum_{e\in E} |w(e)| \\ \\
\d(w,w')=\abs{w-w'} \\ \\
W_0=\os w\in W \svert \abs{w}< \infty \cs \\ \\
W^s_0=W_0 \cap W^s 
\end{array}\end{eqnarray*}
$W$ is a vector space and $W_0$ is a normed subspace of it which, 
under metric $\d$ is isometric to $l_1(\R)$. 
$w,w'\in W$ are isomorphic($w\sim w'$), when there exists a bijection 
$f:V\to V$ such that for every edge $uv\in E$, $w(uv)=w'(f(u),f(v))$.
}

We use standard graphs in constructing unlabeled weighted graphs, 
it has several benefits, in particular with an infinite number of 
isolated vertices we don't need to delete or add a vertex.
\dfn{
Let $\G=W^s/ \sim$ and $\G_0=W_0^s/ \sim$. The elements of $\G$
are named unlabeled graphs. For $w\in W$, we show the equivalence 
class of $w$ by $\wt w$ and define
$$\abs{\wt w}=\abs w$$
}

When no confusion can arise we use the term ``graph" instead of 
``labeled graph" and ``unlabeled graph". 
A graph with all edge weights equal to $0$, simply is denoted by $0$.
We sometimes use common graph theory notions here, converting them to match
our definitions is not difficult, for example ``to delete an edge" means ``to 
change its weight to zero". 
\dfn{
We denote by $S(V)$ the set of all bijections on $V$. 
Suppose $\sigma \in S(V)$ and $e=uv \in E$, let us define
$\sigma(e)=\{\sigma(u),\sigma(v)\}$. 
Also for $w\in W$, define $w^\sigma \in W$ as follows
$$\forall e\in E\quad w^\sigma(\sigma(e))=w(e)$$ 
or equivalently
$$\forall e\in E\quad w^\sigma (e)=w(\sigma^{-1}(e))$$
obviously for $w,w' \in W$
\begin{eqnarray*}\begin{array}{c}
w \sim w' \iff \exists \sigma \in S(V): w'=w^\sigma \\ \\
\wt w=\os w^\sigma\; \vert \;\sigma\in S(V)\cs
\end{array}\end{eqnarray*}
}
\lem{
Given $w,w'\in W$ and $\sigma,\gamma \in S(V)$, the following equations hold
\begin{eqnarray*}\begin{array}{c}
\d(w^\sigma, w'^\sigma)=\d(w,w') \\ \\
{(w^\sigma)}^\gamma=w^{\gamma\circ\sigma}
\end{array}\end{eqnarray*}
}
\dfn{
Given two graphs $G,G'\in \cal G$, we define their distance as follows
$$\d(G,G')=\inf\os \d(w,w')\svert w\in G,w'\in G'\cs $$
Also for $w \in G$ we define
$$\d(G',w)=\d(w,G')=\inf \os \d(w,w')\;\vert \; w' \in G'\cs $$
}
\lem{
For $w_1,w_2\in W^s$
$$\d(\wt w_1,\wt w_2)=\d(w_1,\wt w_2)=\d(\wt w_1,w_2)$$
}
\prf{
\begin{eqnarray*}\begin{array}{c}
\d(w_1^{\sigma_1},w_2^{\sigma_2})=\d(w_1,w_2^{\sigma_1^{-1}\circ\sigma_2})
=\d(w_1^{\sigma_2^{-1}\circ\sigma_1},w_2) \Longrightarrow
\\ \\
\os \d(w_1^{\sigma_1},w_2^{\sigma_2})\svert \sigma_1,\sigma_2\in S(V)\cs =
\os \d(w_1,w_2^{\alpha})\svert \alpha\in S(V))\cs = \\ \\
\os \d(w_1^{\beta},w_2)\svert \beta\in S(V)\cs 
\end{array}\end{eqnarray*}
}
\lem{
For $w_1,w_2,w_3\in W^s$
$$\d(\wt w_1,\wt w_3) \leq \d(\wt w_1,\wt w_2)+\d(\wt w_2,\wt w_3)$$
}
\prf{
$$
\d(\wt w_1,\wt w_2)+\d(\wt w_2,\wt w_3) =
\d(\wt w_1, w_2)+\d(w_2,\wt w_3) \geq
\d(\wt w_1,\wt w_3)
$$
}
\dfn{
Let $w\in W$, $A\subseteq \R$ and $E'\subseteq E$. We define graphs 
$\cut (w,A),\cut(w,E')\in W$ as follows
$$\cut(w,A)(e)=\left\{ \begin{array}{lr}
w(e) & w(e)\in A \\ 0 & {\rm otherwise}
\end{array} \right.$$
this means to delete all edges with weights outside of $A$
$$\cut(w,E')(e)=\left\{ \begin{array}{lr}
w(e) & e\in E' \\ 0 & {\rm otherwise}
\end{array} \right.$$
also for $\epsilon\geq 0$
$$\cut(w,\epsilon)=\cut(w,(-\infty,-\epsilon]\cup[\epsilon,\infty))$$
moreover for $w\in W^s$ we have two more definitions
$$
\cut(\wt w,A)=\wt{\cut(w,A)} \qquad
\text{and} \qquad
\cut(\wt w,\epsilon)=\wt{\cut(w,\epsilon)}
$$
}
\dfn{
The spectrum of a graph is the set of all of it's edge weights
$$\spec(w)=\spec(\wt w)=\os w(e)\svert e \in E\cs $$
}
\lem{
If $w \in W^s_0$ then $\spec(w)$ is a countable and compact set, 
furthermore the only possible limit point of it is $0$.
}
\lem{
If $G,H \in {\cal G}_0$ and $\d(G,H)=0$ then for every $\epsilon \geq 0$, 
$\cut(G,\epsilon)=\cut(H,\epsilon)$.
}
\prf{
Let
$$A=\os |x-y|\svert x\in \spec(G),y\in \spec(H),x\neq y
,(|x|\geq \epsilon\; or\; |y|\geq \epsilon)\cs $$
If $A=\emptyset$ then $\cut(G,\epsilon)=0=\cut(H,\epsilon)$, otherwise we can
define
$$\delta =\min A $$
Clearly $\delta>0$. 
Choose $w\in G$ and $w'\in H$ such that $\abs{w-w'}<\delta$. 
For $e\in E$, $\abs{w(e)-w'(e)}<\delta$, and two cases are possible
\newline
{\bf case i:} $\abs{w(e)},\abs{w'(e)}<\epsilon$, which implies
$$\cut(w,\epsilon)(e)=0=\cut(w',\epsilon)(e)$$
\newline
{\bf case ii:} $w(e)\geq \epsilon$ or $w'(e)\geq\epsilon$, in this case 
$w(e)=w'(e)$ because otherwise $\abs{w(e)-w'(e)}\in A$, 
and hence $\delta<\delta$ which is a contradiction. Thus
$$\cut(w,\epsilon)(e)=w(e)=w'(e)=\cut(w',\epsilon)(e)$$
therefore in each case, $\cut(w,\epsilon)=\cut(w',\epsilon)$ 
and consequently $\cut(G,\epsilon)=\cut(H,\epsilon)$.
}

In the following lemma which is known as Konig infinity lemma, please
forget our notion of a graph, just take it as in ordinary graph theory
texts. 
\lem{
$A_1,A_2,A_3,\cdots$ are nonempty, finite and disjoint sets, 
and $G$ is a graph with $\bigcup_{n=1}^\infty A_n$ as vertex set, 
such that for every $n$, every vertex in $A_{n+1}$ has a neighbour in $A_n$.
$G$ contains a ray $a_1a_2a_3\cdots$ with $a_n\in A_n$. 
(A ray is a sequence of different vertices each of which adjacent to it's successor)
}
\thm{
if $w,w'\in W^s_0$ and $\d(\wt w,\wt w')=0$ then $\wt w=\wt w'$, 
and consequently $d$ is a metric on $\G_0$
}
\prf{
Assume that $w_n=\cut(w,{1\over n})$ and $w'_n=\cut(w',\frac 1 n)$. 
The above lemma implies that $w_n\sim w'_n$. 
Denote by $U_n$ and $U_n'$ the sets of non-isolated vertices of $w_n$ and 
$w'_n$, and by $A_n$ the set of all pairs $(n,f)$ in which $f$ is an 
isomorphism between nonzero parts of $w_n$ and $w'_n$
$$A_n=\{(n,f)\vert f: U_n\leftrightarrow U'_n, \forall u,v \in U_n 
(u\neq v \Rightarrow w_n(uv)=w_n'(f(u),f(v)))\}$$
Since $w_n\sim w'_n$ and $U_n$, $U'_n$ are finite, $A_n$ is nonempty and finite.
Define a graph with vertex set $\bigcup_{n=1}^\infty A_n$ and edge set
$\os \{(n,f),(n+1,g)\}\svert f\subseteq g\cs$. Consider $(n+1,g)\in A_{n+1}$.
Let $f$ be the restriction of $g$ to $U_n$, it is easily seen that 
$(n,f)\in A_{n}$ and $(n,f)$ is a neighbor of $(n+1,g)$ so each vertex 
in $A_{n+1}$ has a neighbour in $A_n$.
Then according to Konig infinity lemma, there is an infinite sequence 
$(1,f_1),(2,f_2),(3,f_3),\cdots$ of vertices such that for each $n$, 
the vertex $(n,f_n)$ is adjacent to the vertex $(n+1,f_{n+1})$, 
i.e $f_1\subseteq f_2\subseteq f_3\subseteq\cdots$. 
We put $f=\bigcup_{n=1}^\infty f_n$.
$f$ is an isomorphism between nonzero parts of $w$ and $w'$.
Since both $w$ and $w'$ have a countable number of isolated vertices, 
$f$ can be extended to an isomorphism between $w$ and $w'$.
}
\section{Completeness of $\G_0$}
In this chapter and the next one we try to find some topological properties of
$\G_0$. 
The main result in this chapter is the completeness of $\G_0$.
\thm{
Let $G,G_1,G_2,\cdots\in {\cal G}_0$. The following are equivalent 
\begin{enumerate}
\item{$G_n\rightarrow G$}
\item{for every $w\in G$ there is a sequence $w_n\in G_n$ such that
$w_n\to w$}
\item{there is a $w\in G$ and a sequence $w_n\in G_n$ such that $w_n\to w$ }
\end{enumerate}
}
\prf{ 
$(1)\Rightarrow(2)$: Let $w\in G$, since 
$$\d(G_n,G)=\d(G_n,w)=\inf \os \d(w',w)\;\vert \; w' \in G_n\cs $$ 
the elements $w_n\in G_n$ exist such that 
$\d(w_n,w)\leq \d(G_n,G)+ \frac{1}{n}$, 
therefore $w_n\to w$. 
$(2)\Rightarrow(3)$ is evident 
\newline
$(3)\Rightarrow(1)$ It follows from the inequality
$$\d(G_n,G)\leq \d(w_n,w)$$
}
\thm{
If the sequence $\{w_n\}\subseteq W^s_0$ be convergent to a graph in $W_0$, 
then the sequence $\{\wt w_n\}$ is convergent in ${\cal G}_0$.
}
\prf{ 
Suppose $w_n \to w$, we shall define the graphs $w'_n$ and $w'$ by means of 
the following equations
$$w_n'(v_iv_j)=\left\{ \begin{array}{lr}
w_n(v_{i/2}v_{j/2}) & i,j\; {\rm are\; even} \\ 
0 & {\rm otherwise}
\end{array} \right.$$
$$w'(v_iv_j)=\left\{ \begin{array}{lr}
w(v_{i/2}v_{j/2}) & i,j\; {\rm are\; even} \\ 
0 & {\rm otherwise}
\end{array} \right.$$
we simply observe that $w'$ is a standard graph and $w_n\sim w_n'$ and $\d(w_n',w')=\d(w_n,w)$ so $w_n' \to w'$, therefore $\wt w_n' \to \wt w'$ and accordingly $\wt w_n \to \wt w'$.
}
\dfn{
Suppose that $w,w'\in W$. 
$w$ is a subgraph of $w'$ ($w\preceq w'$) when 
$$\forall e\in E\quad \Big(w(e)=0\quad\vee\quad w(e)=w'(e)\Big)$$
also we say, $G$ is a subgraph of $G'$ ($G\preceq G'$) where $G,G'\in \cal G$, 
if one of these equivalent statements holds
\begin{eqnarray*}
\exists w\in G \;\exists w'\in G'\quad \Big(w\preceq w'\Big) \\
\forall w\in G \;\exists w'\in G'\quad \Big(w\preceq w'\Big) \\
\forall w'\in G' \;\exists w\in G\quad \Big(w\preceq w'\Big)
\end{eqnarray*}
Clearly if $w,w'\in W_0$ and $w\prec w'$ ،then $\abs w < \abs {w'}$.
Also $\preceq$ is a partial order on $W$.
}
\thm{
$\preceq$ is a partial order on ${\cal G}_0$.
}
\prf{
The transitive and reflexive properties are consequences of the similar 
properties in $W$.
The proof of the antisymmetric property:
If this property fails, then $G\prec G'\prec G$ holds for some 
$G,G'\in\G_0$, implying $\abs G<\abs {G'}<\abs G$, which is impossible.
}
\thm{
If $G,G'\in {\cal G}_0$, then 
\begin{enumerate}
\item{$\d(G,G')\geq \abs{\abs{G}-\abs{G'}}$}
\item{If $G\preceq G'$ then $\d(G,G')=\abs{G'}-\abs{G}$}
\end{enumerate}
}
\prf{
1) For each $w\in G$ and $w'\in G'$ we have
$$\abs{w-w'}\geq \abs{\abs{w}-\abs{w'}}=\abs{\abs{G}-\abs{G'}}$$
so $\d(G,G')\geq \abs{\abs{G}-\abs{G'}}$. 
\newline
2) Choose $w\in G$ and $w'\in G'$ such that
$w\preceq w'$, clearly
$$\abs{w'-w}=\abs{w'}-\abs{w}=\abs{G'}-\abs{G}$$
which in combination with (1) gives the result. 
}

If $G_1\preceq G_2\preceq G_3\preceq\cdots$ is an increasing sequence 
in $\cal G$, then there exists an increasing sequence 
$w_1\preceq w_2\preceq w_3\preceq\cdots$ such that for every $n\in N$, 
$w_n\in G_n$. It is enough to select $w_1$ from $G_1$ then construct 
other terms inductively.
\thm{
Suppose that $w,w_1,w_2,w_3,\cdots\in W_0$, $w_n\to w$ and $\{a_n\}$ 
is a sequence of nonnegative real numbers converging to zero, it follows that 
$\cut(w_n,a_n)\to w$.
}
\prf{
We take $\epsilon>0$ and select $\delta>0$ satisfying 
$\cut(w,[-\delta,\delta])<\epsilon$
and choose $N_1\in N$ such that 
$$\forall n\geq N_1\quad a_n<\delta$$
we also set $A=\os e\in E\svert \abs{w(e)}>\delta \cs$.
Since A is finite, we can select $N_2$ and $N_3$ in such a way that
\begin{eqnarray*}\begin{array}{c}
\forall n\geq N_2\;\forall e\in A\quad \abs{w_n(e)}>\delta \\ \\
\forall n\geq N_3\quad \abs{w_n-w}<\epsilon
\end{array}\end{eqnarray*}
Now for $n\geq\max\{N_1,N_2,N_3\}$
\begin{eqnarray*}\begin{array}{c}
\abs{\cut(w_n,a_n)}\geq \abs{\cut(w_n,A)}\geq\abs{\cut(w,A)}-
\abs{\cut(w-w_n,A)} \\ \\
\geq \abs{w}-\abs{\cut(w,[-\delta,\delta])}-\abs{w-w_n}
\geq \abs{w}-\epsilon-\epsilon=\abs{w}-2\epsilon \\ \\
\Longrightarrow \\ \\
\abs{w-\cut(w_n,a_n)}\leq\abs{w-w_n}+\abs{w_n-\cut(w_n,a_n)} \\ \\
\leq \epsilon+\abs{w_n}-\abs{\cut(w_n,a_n)}\leq\epsilon
+\abs{w}+\epsilon-(\abs{w}-2\epsilon)=4\epsilon
\end{array}\end{eqnarray*}
}

\thm{
In $W_0$ (or $\G_0$), any increasing bounded sequence(with respect to 
$\preceq$) is convergent.
}
\prf{
Suppose that $\{w_n\}$ is an increasing bounded sequence in $W_0$. 
We define the graph $w$ as follows 
$$w(e)=\lim_{n\to \infty} w_n(e)$$ 
The sequence $w_n(e)$ is ultimately constant, so the limit exists. 
It is easily seen that ⃒
$$\abs{w}=\lim_{n\to \infty}\abs{w_n}<\infty\quad\Longrightarrow\quad
w\in W_0$$
Given the fact that، $w_n\preceq w$ for any $n$, we have
$$\lim_{n\to \infty}\d(w,w_n)=\lim_{n\to \infty}(\abs{w}-\abs{w_n})=0$$
therefore $w_n\to w$.
\newline
Now, suppose $\{G_n\}$ is an increasing bounded sequence in $\G_0$. 
Corresponding to this sequence, there is an increasing sequence  
$\{w_n\}\subseteq W_0$ such that $w_n\in G_n$.
The convergence of $\{G_n\}$ is a result of the convergence of 
$\{w_n\}$.
}
\dfn{
A graph in which no two nonzero edges are adjacent is called
a jointless one. 
We denote by $\wh\G$ and $\wh\G_0$ the sets of jointless graphs in $\G$ 
and $\G_0$.
Given a graph $G\in \G_0$, there is a unique member of $\wh\G_0$ which has
the same edge weights as $G$ (with same multiplicity), we denote it by $\wh G$,
see example 4.
}
\dfn{
Suppose $G\in \wh\G$ and $a\in \spec(G)$. The unique graph obtained 
from $G$ by deleting an edge with weight $a$ is denoted by $\del(G,a)$. 
}
\dfn{
Suppose $\wt{w}=G\in \G_0$. We define $\code(w)=\code(G)=f$
where $f\in l_1(\R)$ is constructed by induction as follows: set $G_0=G$ and 
define
$$
f(n)=\left\{ \begin{array}{lr}
\max\;\spec(G_{n-1}) & n\; {\rm is\; odd} \\ 
\min\;\spec(G_{n-1}) & n\; {\rm is\; even}
\end{array} \right.
\qquad G_n=\del(G_{n-1},f(n))
$$
}\label{defcode}
\xmpl{
Suppose that $G$ is the jointless graph that has an edge with weight 
$\frac{1}{2^n}$ for each $n\geq 1$ and two edges of weight $-1$ , then
$$\code(G)=(\frac{1}{2},-1,\frac{1}{4},-1,\frac{1}{8},0,
\frac{1}{16},0,\cdots)$$
}
\lem{
Suppose $G,H\in \wh\G_0$, $x_1=\max\spec(G)$ and $x_2=\max\spec(H)$ 
(or $x_1=\min\spec(G)$ and $x_2=\min\spec(H)$ ) and
$G'=\del(G,x_1)$ and $H'=\del(H,x_2)$, then
$$\d(G,H)=\d(G',H')+\abs{x_1-x_2}$$
}
\prf{
It is clear that
$\d(G,H)\leq \d(G',H')+\abs{x_1-x_2}$. To show the
inverse, we take
$w_1\in G$ and $w_2\in H$. Let $x_1=w_1(e_1)$, $x_2=w_2(e_2)$, 
$y_1=w_1(e_2)$ and $y_2=w_2(e_1)$ and define the graph $w$ as follows
$$
w(e)=\left\{ \begin{array}{lr}
x_1 & e=e_2 \\ 
y_1 & e=e_1 \\
w_1(e) & {\rm otherwise}
\end{array} \right.
$$
according to the assumption, $y_1\leq x_1$ and $y_2\leq x_2$ (or $x_1\leq y_1$ and 
$x_2\leq y_2$ ). Using these relations one can easily show that
\begin{eqnarray*}\begin{array}{c}
\d(w_1,w_2)-\d(w,w_2)=\abs{x_1-y_2}+\abs{x_2-y_1}-\abs{x_1-x_2}-\abs{y_2-y_1}
\geq 0 \\ \\
\Longrightarrow\quad
\d(w_1,w_2)\geq \d(w,w_2)\geq \abs{x_1-x_2}+\d(G',H')
\end{array}\end{eqnarray*}
}
\thm{
For every $G,H\in \G_0$ we have 
$$\d(G,H)=\abs{\code(G)-\code(H)}$$
}
\prf{
Let $G_n$ be the sequence related to graph $G$ in definition 10, 
and relate a similar sequence $H_n$ to $H$.
We set $\code(G)=g$ and $\code(H)=h$. Applying the preceding 
lemma $n$ times, we obtain
$$\d(G,H)=\sum_{i=1}^n \abs{g(n)-h(n)}+\d(G_n,H_n)$$
Since the sequences $\abs{G_n}$ and $\abs{H_n}$ are convergent to $0$, 
$\d(G_n,H_n)\to 0$ and consequently
$$\d(G,H)=\sum_{i=1}^\infty \abs{g(n)-h(n)}=\abs{\code(G)-\code(H)}$$
}
\lem{
$\widehat\G_0$ is a complete subspace of $\G_0$
}
\prf{
Set $B=\os \code(G) \svert G\in \wh\G_0 \cs$ 
it is easily seen that 
\begin{eqnarray*}\begin{array}{c}
B=\os f\in l_1 \svert \forall n\; \Big(f(2n-1)\geq 0,f(2n)\leq 0,
\\ \\ \qquad\quad f(2n-1)\geq f(2n+1),f(2n)\leq f(2n+2) \Big)\cs
\end{array}\end{eqnarray*}
and B is a closed subset and consequently a complete subset of $l_1(\R)$.
According to the theorem 8, $\code: \wh\G_0 \to B$ is an onto isometry, 
so $\wh\G_0$ is also complete. 
}

By a correspondance between two graphs $G_1,G_2\in \G$ we mean a choice of two 
members $w_1\in G_1$ and $w_2\in G_2$.
\thm{
$\G_0$ is a complete metric space 
}
\prf{
Suppose that $G_n$ is a Cauchy sequence. Since 
$\d(\wh G_n,\wh G_m)\leq \d(G_n,G_m)$, 
so the sequence $\{\wh G_n\}$ is also Cauchy. Therefore according to the 
preceding lemma, the sequence $\{\wh G_n\}$ is convergent to a jointless 
graph $F$. Set $A=\spec(F)$ and  
let $w_n \in G_n$.
Suppose $w'_n$ is obtained from $w_n$ by 
rounding the weight of each edge to the closest number in $A$ (if the weight 
of an edge has the least difference with two numbers in A, we choose one 
of them arbitrarily). 
Set $H_n=\wt w'_n$,
we claim that $\d(H_n,G_n)\to 0$. In fact, it is easily seen that 
$\d(G_n,H_n)\leq \d(w_n,w'_n)\leq \d(\wh G_n,F)$, which proves the claim.
Therefore, it is enough to show $\{H_n\}$'s convergence instead of $\{G_n\}$'s.
First let us show that for all $\epsilon>0$ there exists a natural number $M$ 
such that 
\begin{eqnarray}\label{e1}
\forall m\geq M\quad \cut(H_m,\epsilon)=\cut(H_M,\epsilon)
\end{eqnarray}
To prove that we set $A_\epsilon=\os x\in A \svert \abs{x}\geq \epsilon \cs$ and 
$\delta=\min \os \d(x,A\backslash \{ x \}) \svert x\in A_\epsilon \cs$ .
Since $A$ does not have a nonzero limit point and $A_\epsilon$ is finite, 
so $\delta>0$. 
The relation $\d(H_n,G_n)\to 0$ shows that $\{H_n\}$ is a Cauchy sequence 
so there is a $M\in \mathbb N$ such that $\d(H_M,H_m)<\delta$ for every 
$m\geq M$. 
Now, if $\cut(H_m,\epsilon)\neq \cut(H_M,\epsilon)$ for one $m\geq M$, 
in every correspondence between $H_m$ and $H_M$ we get an edge which has 
two different weights in the two graphs, one of which from $A_\epsilon$, 
and the other from $A$, and hence
$\d(H_M,H_m)\geq \delta$, which is a contradiction. 
Choose a strictly increasing sequence $\{M_n\}$ such that for each $n$, 
$\epsilon=\frac{1}{n}$ and $M=M_n$ satisfy the equation (\ref{e1}). It is 
evident that for each $n$, 
$\cut(H_{M_n},\frac{1}{n})\preceq\cut(H_{M_{n+1}},\frac{1}{n+1})$.
Then, since $\abs{\cut(H_{M_n},\frac{1}{n})}\leq \abs{H_{M_n}}$, so the 
sequence $\{\cut(H_{M_n},\frac{1}{n})\}$ is bounded and consequently 
converges to a graph $G$. We have
\begin{eqnarray}\label{e2}\begin{array}{c}
\d(\wh H_n,\wh G_n)\leq \d(H_n,G_n) \Longrightarrow \d(\wh H_n,\wh G_n)\to 0\\ \\
\Longrightarrow \wh H_n \to F \Longrightarrow \wh H_{M_n} \to F
\Longrightarrow\\ \\ 
\cut(\wh H_{M_n},\frac{1}{n})\to F \Longrightarrow 
\d(\wh H_{M_n},\cut(\wh H_{M_n},\frac{1}{n}))\to 0
\end{array}\end{eqnarray}
Also
\begin{eqnarray}\label{e3}\begin{array}{c}
\d(H_{M_n},\cut(H_{M_n},\frac{1}{n}))=\abs{H_{M_n}}-\abs{\cut(H_{M_n},\frac{1}{n})}
\\ \\=\abs{\wh H_{M_n}}-\abs{\cut(\wh H_{M_n},\frac{1}{n})}
=\d(\wh H_{M_n},\cut(\wh H_{M_n},\frac{1}{n}))
\end{array}\end{eqnarray}
(\ref{e2}) and (\ref{e3}) conclude that
$$\d(H_{M_n},\cut(H_{M_n},\frac{1}{n}))\to 0$$
Therefore $H_{M_n}\to G$, so $H_n$ has a convergent subsequence and
consequently, it is convergent itself.
}
\section{Examining the space for some other common topological properties}
Besides the completeness of $\G_0$ there are some other important
properties of the space, we discuss a few of them here. Note
that, in this chapter by a finite graph we mean one that has only
finitely many non-isolated vertices. Also in a metric space
$\mathcal{M}$ we denote the ball with center $x$ and radius $r$ by
$B_{\mathcal{M}}(x,r)$. \thm{ $\G_0$ is separable. } \prf{ The set
of all finite graphs with rational edge weights is a dense subset of
$\G_0$. }
By a finite graph we mean one that has only finitely many non-isolated vertices.
Also in a metric space $\mathcal{M}$ we denote the ball with
center $x$ and radius $r$ by $B_{\mathcal{M}}(x,r)$.
\thm{
$\G_0$ is separable.
}
\prf{
The set of all finite graphs with rational edge weights is a dense subset 
of $\G_0$.
}
\thm{
$\G_0$ is not locally compact. 
}
\prf{
Suppose the contrary, so there is an open neighbourhood $B_{\G_0}(0,r)$  
such that $\overline{B_{\G_0}(0,r)}$ is compact. 
Consider the sequence $\{ G_n \}$ in which $G_n$ is the jointless graph with 
$n$ edges of weight $\frac{r}{n}$.
This sequence must have a convergent subsequence.
The tiny edges of this subsequence say that it converges to $0$. On the 
other hand the norm of its members are always equal to $r$ implying that the 
norm of the limit graph must be $r$, which is a contradiction.
}
\lem{
$W_0^s$ is path connected.
}
\prf{
Take the standard graph $w$. The function $f:[0,1]\to W_0^s$, $f(t)=tw$ is
a path between 0 and $w$, so every point is connected to 0 via a path.
}
\thm{
$\G_0$ is path connected.
}
\prf{
The onto function $w\to \wt{w}$ 
from $W_0^s$ to $\G_0$ is continuous and hence takes path connected
to path connected.
}
\thm{
$\G_0$ is locally path connected.
}
\prf{
First we show that every ball in $W_0^s$ is path connected.
Let $w\in W_0^s$, $r > 0$ and $w'\in B_{W_0^s}(w,r)$. 
Choose a finite graph $w''$ from $B_{W_0^s}(w,r)$ (finite graphs
are dense in $W_0^s$).
$B_{W_0}(w,r)$ is convex and the graphs $w'$ and $w''$ have 
infinitely many isolated vertices in common, so we can define 
the function $f:[0,1]\to B_{W_0^s}(w,r)$, $f(t)=tw'+(1-t)w''$
which is a path between $w'$ and $w''$. So every member of 
$B_{W_0^s}(w,r)$ is connected to $w''$ via a path, therefore
the ball is path connected.
Now according to the facts that for each $w\in W_0$ we have 
$B_{\G_0}(\wt w,r)=\wt{B_{W_0^s}(w,r)}$ (which is not difficult to prove)
and that the function $w\to \wt{w}$ 
from $W_0$ to $\G_0$ is continuous, every ball in $\G_0$ is path 
connected and so $\G_0$ is locally path connected.
}

\bigskip\noindent {\bf Acknowledgments.}
We are grateful to professor Sayyad Ebadollah Mahmoodian, for his valuable 
advices and comments during this work. We also thank all friends
that directly and indirectly gave us hands in the preparation of this paper.

\end{document}